\documentclass[11pt]{article}
\usepackage{cite}
\usepackage{mathrsfs}
\usepackage{amsfonts}
\usepackage{amsmath}
\usepackage{amsfonts,amssymb,color}
\usepackage{dsfont}
\usepackage{curves}
\usepackage{mathrsfs}
\usepackage{pifont}
\usepackage{amssymb}
\usepackage{latexsym,amsmath,amssymb,amsfonts,epsfig,graphicx,cite,psfrag}
\usepackage{eepic,color,colordvi,amscd}

\newtheorem{theorem}{Theorem}[section]
\newtheorem{proposition}{Proposition}[section]

\newtheorem{lemma}{Lemma}[section]
\newtheorem{corollary}{Corollary}[section]

\newtheorem{claim}{Claim}[section]
\newtheorem{conjecture}{Conjecture}[section]

\newcommand{\qed}{\hfill\rule{0.5em}{0.809em}}

\def\emptyset{\mbox{{\rm \O}}}

\renewcommand{\baselinestretch}{1.1}
\def\qed{\hfill \rule{4pt}{7pt}}

\def\pf{\noindent {\it Proof. }}

\begin{document}

\title{On coloring of graphs of girth $2\l+1$ without longer
odd holes\thanks{Partially supported by NSFC projects 11931006. This paper is published in SCIENCE CHINA Mathematics (in Chinese) Doi: 10.1007/s11425-021-1949-8}}
 \author{Di Wu$^{1,}$\footnote{Email: 1975335772@qq.com},  \;\; Baogang  Xu$^{1,}$\footnote{Email: baogxu@njnu.edu.cn OR baogxu@hotmail.com.}\;\; and \;\; Yian Xu$^{2,}$\footnote{Email: yian$\_$xu@seu.edu.cn}\\\\
\small $^1$Institute of Mathematics, School of Mathematical Sciences\\
\small Nanjing Normal University, 1 Wenyuan Road,  Nanjing, 210023,  China\\
\small $^2$School of Mathematics, Southeast University, 2 SEU Road, Nanjing, 211189, China}
\date{}

\maketitle


\begin{abstract}
A hole is an induced cycle of length at least 4. Let $\l\ge 2$ be a positive integer, let ${\cal G}_l$ denote the family of graphs which have girth $2\l+1$ and have no holes of odd length at least $2\l+3$, and let $G\in {\cal G}_{\l}$.  For a vertex $u\in V(G)$ and a nonempty set $S\subseteq V(G)$,  let $d(u, S)=\min\{d(u, v):v\in S\}$, and let $L_i(S)=\{u\in V(G) \mbox{ and } d(u, S)=i\}$ for any  integer $i\ge 0$. We show that if $G[S]$ is connected and $G[L_i(S)]$ is bipartite for each $i\in\{1, \ldots, \lfloor{\l\over 2}\rfloor\}$, then $G[L_i(S)]$ is bipartite for each $i>0$, and consequently $\chi(G)\le 4$, where $G[S]$  denotes the subgraph induced by $S$.
Let $\theta^-$ be the graph obtained from the Petersen graph by deleting three vertices which induce a path, let $\theta^+$ be the graph obtained from the Petersen graph by deleting two adjacent vertices, and let $\theta$ be the graph obtained from $\theta^+$ by removing an edge incident with two vertices of degree 3. For a graph $G\in{\cal G}_2$, we show that if $G$ is 3-connected and has no  unstable 3-cutset then $G$ must induce either $\theta$ or $\theta^-$ but does not induce $\theta^+$. As corollaries, $\chi(G)\le 3$ for every graph $G$ of ${\cal G}_2$ that induces neither $\theta$ nor $\theta^-$, and minimal non-3-colorable graphs of ${\cal G}_2$ induce no $\theta^+$.

\begin{flushleft}
\end{flushleft}

\end{abstract}

\newpage

\section{Introduction}

Let $G$ be a graph, and let $u$ and $v$ be two vertices of $G$. We simply write $u\sim v$ if $uv\in E(G)$, and write $u\not\sim v$ if $uv\not\in E(G)$. Let $S$ be a subset of $V(G)$. We use $G[S]$ to denote the subgraph of $G$ induced by $S$, call $S$ a {\em clique} if $G[S]$ is complete, call $S$ a {\em stable set} if $G[S]$ has no edges, and say that $S$ is {\em unstable} if it is not stable. For two induced subgraphs $F$ and $H$ of $G$, in the case that it makes no confusions, we simply write $V(F)$ as $F$, and  write $G[V(F)\cup V(H)]$ as  $F\cup H$.

Let $S$ and $T$ be two subsets of $V(G)$, and let $x$ and $y$ be two vertices of $G$.  We use $N_S(x)$ to denote the neighbors of $x$ in $S$, and define $N_S(T)=\cup_{x\in T} N_S(x)$ (if $S=V(G)$ then we omit the subindex and simply write $N(x)$ or $N(T)$). If $S\cap T=\emptyset$, we use $E(S, T)$ and $e(S, T)$ to denote the set and number of edges, respectively, with one end in $S$ and the other in $T$. An $xy$-path refers to a path between $x$ to $y$, and an $(S, T)$-path means a path $P$ with  $|S\cap P|=|T\cap P|=1$.   A vertex of degree $k$ is  called a $k$-{\em vertex}, a cycle of length $k$ is  called a $k$-{\em cycle}, and a cutset of size $k$ is called a $k$-cutset.  For a path $P$, we use $\l(P)$ and $P^*$ to denote the length and the set of internal vertices of $P$, respectively. An {\em odd} (resp. {\em even}) path refers to a path of odd (resp. even) length.

Let $H$  be a graph. We say that $G$ induces $H$ if $G$ has an induced subgraph isomorphic to $H$, and say that $G$ is $H$-{\em free} if $G$ does not induce $H$.
A {\em hole}  is an induced cycle of length at least 4, and an {\em odd} (resp. {\em even})  hole is a hole of  odd length (resp. even length).

Let $k$ be a positive integer, and let $[k]=\{1, 2, \ldots, k\}$. A $k$-{\em coloring} of $G$ is a mapping $c: V(G)\mapsto [k]$ such that $c(u)\neq c(v)$ whenever $u\sim v$ in $G$.  The {\em chromatic number} $\chi(G)$ of $G$ is the minimum integer $k$ such that $G$ admits a $k$-coloring. It is certain that $\chi(G)\ge \omega(G)$. But the difference $\chi(G) - \omega(G)$ may be arbitrarily large as there are triangle-free graphs with arbitrary large chromatic number (see \cite{BD1947, AAZ1949,JM1955}). Where $\omega(G)$ denotes the {\em clique number} which is the maximum size of cliques of $G$. Erd\H{o}s \cite{PE1959} showed that for every pair of positive integers $k\ge 2$ and $l\ge 3$ there exists a graph $G$ with $\chi(G)=k$ whose shortest cycle has length at least $l$.

For a given family ${\cal G}$ of graphs, if there is a function $f$ such that $\chi(G)\le f(\omega(G))$ for each graph $G\in {\cal G}$, then the family ${\cal G}$ is said to be $\chi$-{\em bounded},  and the function $f$ is called a {\em binding function} of ${\cal G}$  (Gy\'{a}rf\'{a}s, \cite{gyarfas1}).

The $\chi$-boundedness of graphs  with some particular forbidden holes were studied extensively. Addario-Berry, Chudnovsky, Havet, Reed and  Seymour \cite{ACHRS08}, and Chudnovsky and  Seymour \cite{MCPS2019} fixed a gap of \cite{ACHRS08}, confirmed a conjecture of Reed and proved that every even hole free graph has a vertex whose neighbors is the union of two cliques. As a direct consequence, every even hole free graph $G$ has $\chi(G)\le 2\omega(G)-1$. One may find more results and problems about even hole free graphs in \cite{KV2010}.

Let  $\l$ be a positive integer. Confirming three conjectures of Gy\'{a}rf\'{a}s \cite{gyarfas1987}, Scott and Seymour \cite{SS2015} proved that
$\chi(G)\le 2^{2^{\omega(G)+2}}$  for odd  hole free graphs, Chudnovsky, Scott and  Seymour \cite{MCSS2017} proved that (holes of length at least $l$)-free graphs are $\chi$-bounded, and Chudnovsky, Scott, Seymour and Spirkl \cite{MSSS2017} further showed that (odd holes of length at least $l$)-free graphs are $\chi$-bounded. Let ${\cal H}$ be the family of graphs with neither triangle nor hole of length 0 modulo 3. Answering a question of Kalai and Meshulam, Bonamy, Charbit and Thomass\'{e} \cite{BCST20} proved that there is a constant $c$ such that $\chi(G)\le c$ for each graph $G\in {\cal H}$. The question that whether $\chi(G)\le 3$ for each $G\in {\cal H}$ remains open. There are also quite a lot of results concerning the structure and chromatic number of graphs inducing no paths on $\l$ vertices, one may see \cite{BCMZ18, CHPS20, KCHM19, CMSZ2017, CKMM20, SGSH19, SHTK19, TKFM19} for some most recent results. Interested readers are referred to \cite{ISBR19} and \cite{SS2020} for more information on $\chi$-bounded problems.

In this paper, we study a family  ${\cal G}_{\l}$ of graphs that have girth $2\l+1$ and have no odd holes of length at least $2\l+3$, where the {\em girth} of a graph is the length of a shortest cycle in it.  Robertson conjectured (see \cite{NPRZ2011}) that the only 3-connected, internally 4-connected  graph in ${\cal G}_2$ is the Petersen graph. Plummer and Zha \cite{MPXZ} presented counterexamples to Robertson's conjecture, and posed the following conjecture.

\begin{conjecture}\label{conj-1} {\em (\cite{MPXZ})}
$\chi(G)\le 3$ for each graph $G$ of ${\cal G}_2$.
\end{conjecture}

In \cite{XYZ2017}, Xu, Yu and Zha  proved that

\renewcommand{\baselinestretch}{1}
\begin{theorem}\label{theo-1}{\em (\cite{XYZ2017})}
Let $G$ be a graph in ${\cal G}_2$, and let $u$ be a vertex of $G$. Then,
the set of vertices of the same distance to $u$ induces a bipartite subgraph, and consequently
$\chi(G)\le 4$.
\end{theorem}\renewcommand{\baselinestretch}{1.1}

Let $S$ be a set of vertices, and let $x$ be a vertex. We define $d(x, S)=\min\{d(x, y): y\in S\}$, and define $L_i(S)=\{x\in V(G), d(x, S)=i\}$ for any integer $i\ge 0$ (we simply write $L_i(u)$ if $S$ consists of a single vertex $u$).

In this paper, we generalize Theorem~\ref{theo-1} to $\cup_{\l\ge 2} {\cal G}_{\l}$.

\renewcommand{\baselinestretch}{1}
\begin{theorem}\label{lem-subgraph}
Let $\l\ge 2$ be an integer, let $G$ be  a graph in ${\cal G}_{\l}$, and let $S\subset V(G)$ with $G[S]$ connected. Then, $G-e\in{\cal G}_{\l}$ for each edge $e$ not contained in any $(2\l+1)$-cycle, $G[L_i(S)]$ is bipartite for every $i>0$ if $G[L_i(S)]$ is bipartite for each $i\in [\lfloor{\l\over 2}\rfloor]$, and $\chi(G)\le 4$.
\end{theorem}\renewcommand{\baselinestretch}{1.1}

By Theorem~\ref{theo-1}, $\chi(G)\le 4$ for every $G\in{\cal G}_2$. If Conjecture~\ref{conj-1}  is not true, there would be a graph in ${\cal G}_2$ with chromatic number 4. We define
$${\cal G}_0=\{G\in {\cal G}_2 \mbox{ such that $\chi(G)=4$ but each proper subgraph of $G$ is 3-colorable} \}.$$

\renewcommand{\baselinestretch}{1}
\begin{theorem}\label{4-critical}
Let $G$ be a graph in ${\cal G}_0$, and let $u\in V(G)$. Then, every edge is contained in some $5$-hole, and $\delta(G)\ge 4$. Furthermore,
\begin{itemize}
\item [$(a)$] $G$ is $3$-connected,

\item [$(b)$] $\{u, u_1, u_2, u_3\}$ is not a cutset for any three vertices
       $u_1, u_2, u_3\in N(u)$, and

\item [$(c)$] every $3$-cutset of $G$ is stable.

\end{itemize}
Where $\delta(G)$ denotes the minimum degree of $G$.
\end{theorem}\renewcommand{\baselinestretch}{1.1}

\renewcommand{\baselinestretch}{1}
To support Conjecture~\ref{conj-1}, we show that
\begin{theorem}\label{coro-4-colorable}
Graphs in ${\cal G}_2$ which induce no two $5$-cycles sharing edges are $3$-colorable.
\end{theorem}\renewcommand{\baselinestretch}{1.1}

It is easy to see that Theorem~\ref{coro-4-colorable} is a direct consequence of Theorem~\ref{4-critical} and the following

\begin{theorem}\label{coro-4-colorable-1}
Let $G$ be a $3$-connected graph of ${\cal G}_2$ that has a $5$-cycle. If every $3$-cutset of $G$ is stable, then $G$ induces two $5$-cycles sharing edges.
\end{theorem}

Assuming Theorems\ref{4-critical} and \ref{coro-4-colorable-1}, we can now prove Theorem~\ref{coro-4-colorable}.

\noindent{\em Proof of} Theorem~\ref{coro-4-colorable}: If it is not true, we choose $G$ to be a smallest counterexample. It is certain that $G$ must have 5-cycles. Let $C=x_0x_1x_2x_3x_4x_0$ be a 5-cycle of $G$. By Theorem~\ref{4-critical}, $G$ is 3-connected, $\{x_i, x_{i+1}, x\}$ is not a cutset for each $i\in \{0, 1, 2, 3, 4\}$ and each $x\in V(G)$. This contradiction to  Theorem~\ref{coro-4-colorable-1} shows that $\chi(G)\le 3$. \qed

\medskip

Following from Theorem~\ref{coro-4-colorable}, to confirm Conjecture~\ref{conj-1}, we need only to consider graphs of ${\cal G}_2$ that induce two 5-cycles sharing edges.
We use $\theta^-$ to denote the graph obtained from the Petersen graph by removing three vertices which induce a path, use $\theta^+$ to denote the graph obtained from the Petersen graph by removing two adjacent vertices, and use $\theta$ to denote the graph obtained from $\theta^+$ by removing an edge incident with two 3-vertices.
It is easy to check that if a graph $G\in {\cal G}_2$ induces two 5-cycles sharing edges, then it must induce either $\theta$ or $\theta^-$. So, Conjecture~\ref{conj-1} is equivalent to the following

\begin{conjecture}\label{conj-SWXX}
Graphs in ${\cal G}_0$ induce neither $\theta$ nor $\theta^-$.
\end{conjecture}

\begin{figure}[htbp]\label{fig-3}
  \begin{center}
    \includegraphics[width=8cm]{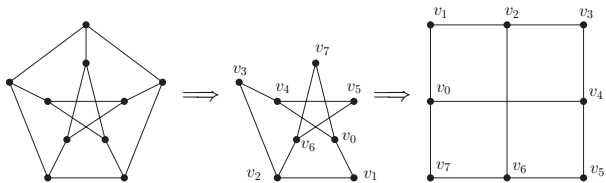}
  \end{center}
  \vskip -15pt
\caption{Getting $\theta^+$ from the Petersen graph}
\end{figure}

Note that $\theta^-$ can be obtained from $\theta^+$ by removing a 2-vertex, and $\theta$ can be obtained from $\theta^+$ by removing an edge joining two 3-vertices. As a support to Conjecture~\ref{conj-SWXX} (also a further support to  Conjecture~\ref{conj-1}), we prove the following

\begin{theorem}\label{4-critical1}
Let $G$ be a $3$-connected graph in ${\cal G}_2$ of which every $3$-cutset is stable. If $G$ is not the Petersen graph, then it is $\theta^+-$free.
\end{theorem}\renewcommand{\baselinestretch}{1.1}

As a direct consequence of Theorems~\ref{4-critical} and \ref{4-critical1}, we have
\begin{corollary}
Graphs in ${\cal G}_0$ are $\theta^+-$free.
\end{corollary}

Let $P$ be a path, and let $x$ and $y$ be two vertices on $P$. We use $P[x, y]$ to denote the segment of $P$ between $x$ and $y$. For a cycle $C=u_0u_1\ldots u_hu_0$ and two distinct integers $0\le i, j\le h$, we use $C[u_i, u_j]$ to denote the path $u_iu_{i+1}\ldots u_j$ (here and somewhere else in the paper, the summation of subindices is always taken modulo the length of the cycle). Let $x$ and $y$ be two vertices of $G$, and $H$ a connected subgraph of $G-\{x, y\}$ such that both $x$ and $y$ have neighbors in $H$. An $xy$-$H$-path is an $xy$-path with all internal vertices in $H$.

In Sections 2 and 3, we prove, respectively,  Theorem~\ref{lem-subgraph} and \ref{4-critical}. Then, we prove Theorem~\ref{coro-4-colorable-1} in Section 4, and prove Theorem~\ref{4-critical1} in Section 5.

\section{Graphs in ${\cal G}_{\l}$ with $\l\ge 2$}

In this section, we focus on the structure on graphs in ${\cal G}_{\l}$ with $\l\ge 2$, and prove Theorem~\ref{lem-subgraph}.

Let $e$ be an edge of $G$ with $G-e\not\in {\cal G}_{\l}$. It is certain that $G-e$ has girth at least $2\l+1$, and so $G-e$ has an odd hole of length larger than $2\l+1$, say $C=u_0u_1\ldots u_{2h}u_0$, and  $e$ must be the unique chord of $C$ in $G$. We may suppose that  $e=u_0u_j$ by symmetry. Now, the parity of $C$ indicates that either $u_0u_1\ldots u_ju_0$ or $u_ju_{j+1}\ldots u_{2h}u_0u_j$ is an odd hole of length $2\l+1$. Therefore, $e$ is contained in  some $(2\l+1)$-cycle.

For simplicity, we write $L_i(S)$ as $L_i$ for all $i\ge 0$, and suppose that $G[L_i]$ is bipartite for each $i\in [\lfloor{\l\over 2}\rfloor]$.
For integers $0\le i<j$ and vertices $u\in L_i$ and $v\in L_j$, we call $u$ an $S$-{\em predecessor} of $v$ (resp.  $v$ an $S$-{\em successor} of $u$) if there exists a $uv$-path of length $j-i$. Let $W$ be a set of vertices. A vertex $x$ is called an $S$-predecessor of $W$ if $x$ is an $S$-predecessor of some vertex of $W$.

To prove the second statement, we suppose that, for some $h\ge \lfloor{\l\over 2}\rfloor$, $G[L_i]$ is bipartite for each $i\in [h]$, but the graph $H=G[L_{h+1}]$ is not bipartite. Then, $H$ must have a $(2\l+1)$-cycle, say $C=u_0u_1\ldots u_{2\l}u_0$.

Let $q=\lfloor{\l\over 2}\rfloor$. Since $G$ has girth $2\l+1$, we see that $u_i$ and $u_j$ ($0\le i< j\le 2\l$) have no common $S$-predecessor in $\cup_{r=h-q+1}^{h} L_r$, otherwise there would be a cycle of length at most $2\l$.

Let $F=G[\cup_{i=0}^{h-q} L_i]$, let $W$ be a minimal set such that every vertex of $C$ has an $S$-predecessor in $W\cap L_r$ for each $r\in \{h-q+1, \ldots, h\}$. Since $G$ has girth $2\l+1$, we see that $G[W]$ consists of  $2\l+1$ induced $(q-1)$-paths.  For each  $0\le i\le 2\l$, let $w_i$ be the $S$-predecessor of $u_i$ in $W\cap L_{h-q+1}$,  let $v_i$ be an $S$-predecessor of $w_i$ in $L_{h-q}$, and let $P_i$ be the  $u_iv_i$-path of length $q+1$ with $P^*_i\subseteq W$.

Suppose that $\l$ is even. If $v_0\sim w_{\l}$, then $v_0w_{\l}(P_{\l}-v_{\l})C[u_{\l}, u_{0}]P_0$ is a $2\l+3$-hole. So, we have that $v_0\not\sim w_{\l}$. By symmetry, we see that $w_0\not\sim v_{\l}$. If $v_0\sim v_{\l}$, then $P_{0}C[u_0, u_{\l}]P_{\l}v_{\l}v_0$ is a $2\l+3$-hole. This shows that $v_0\not\sim v_{\l}$. Let $Q_0$ be a shortest $v_0v_{\l}$-$F$-path. Then, we have $Q_{0}P_{0}C[u_{0}, u_{\l}]P_{\l}$ or $Q_{0}P_{\l}C[u_{\l}, u_{0}]P_{0}$ is an odd hole of length greater than $2\l+1$, a contradiction. So, we have that $\l$ must be odd, and is at least 3. If $v_0\sim w_{\l-1}$, then $v_0P_0u_0u_{1}\ldots u_{\l-1}P_{\l-1}v_0$ is a $2\l$-hole, a contradiction. Therefore, $v_0\not\sim w_{\l-1}$, and $w_0\not\sim v_{\l-1}$ by symmetry. If $v_0\not\sim v_{\l-1}$, let $Q_1$ be a shortest $v_0v_{\l-1}$-$F$-path, then $Q_{1}P_{0}C[u_{0}, u_{\l-1}]P_{\l-1}$ or $Q_{1}P_{\l-1}C[u_{\l-1}, u_{0}]P_{0}$ is an odd hole of length greater than $2\l+1$. So, we have that $v_0\sim v_{\l-1}$. With the same argument, we can show that for any $i\in \{0, 1, 2\l\}$, $v_i\sim v_{i+\l-1}$ and $v_i\sim v_{i+\l+2}$, which imply that $G[\{v_0, v_1, \ldots, v_{2\l}\}]$ has cycles, contradicting the fact that $G[L_{h-q}]$ is bipartite and $G$ has girth $2\l+1$.
This contradiction shows that $G[L_i]$ is bipartite for every $i>0$.

Let $u$ be an arbitrary vertex of $G$. By taking $S=\{u\}$, we see that $L_i(u)$ is a stable set for each $i\in [\l-1]$, and so $L_i(u)$ is bipartite for each integer $i\ge 0$. A 4-coloring of $G$ can be obtained by properly coloring all vertices of $L_i(u)$ with colors 1 and 2 if $i$ is even, and coloring all vertices of $L_i(u)$ with colors 3 and 4 if $i$ is odd.  \qed

\section{Proof of Theorem~\ref{4-critical}}

We prove Theorem~\ref{4-critical} in this section.
Recall that ${\cal G}_0$ denotes the set of 4-chromatic critical graphs in ${\cal G}_2$.

Let $G$ be a graph, let $x$ and $y$ be two vertices of $G$, and let $f$ be a $k$-coloring of $G$. For two distinct integers $i, j\in [k]$, an $xy$-path is called an $(i, j)$-$xy$-{\em path} if all its vertices are colored $i$ and $j$ alternatively. Note that an $(i, j)$-$xy$-path has even length if $f(x)=f(y)$, and has odd length otherwise. This will be used frequently.

In the following, we always suppose that $G$ is a graph in  ${\cal G}_0$. Then, $\delta(G)\ge 3$, and $G$ has no cutset which is a clique.
For simplicity, we call an odd hole of length at least 7 as a {\em big-odd-hole}.

Suppose that $G$ has an edge, say $uv$, that is not contained in any $C_5$. By Theorem~\ref{lem-subgraph}, $G-uv$ is still in ${\cal G}_2$, and thus $\chi(G-uv)=3$ by the choice of $G$. Let $(V_1, V_2, V_3)$ be a partition of $V(G-uv)$ into three stable sets. Since $\chi(G)=4$, we may suppose, by symmetry, that $u, v\in V_1$. Now, $G-V_3$ has 5-cycles, and every odd cycle of $G-V_3$ must contain $uv$, a contradiction. So, every edge of $G$ is in some $C_5$.

If $\delta(G)=3$, choose $u\in V(G)$ with $N(u)=\{u_1, u_2, u_3\}$, and let $\phi$ be a 3-coloring of $G-u$.  Since  $\chi(G)=4$, we may suppose that $\phi(u_i)=i$. Now, $G-u$ has  $(i, j)$-$u_iu_j$-paths for each pair of distinct integers $i, j\in[3]$. Let $P_{i, j}$ be a shortest $(i, j)$-$u_iu_j$-path. Since $G$ has no big-odd-hole and since $u_1, u_2, u_3\in N(u)$, we see that $\l(P_{1, 2})=\l(P_{2, 3})=\l(P_{1, 3})=3$. But now, $P_{1, 2}\cup P_{2, 3}\cup P_{1, 3}$ is a big-odd-hole. Therefore, $\delta(G)\ge 4$.

Now we show that $G$ is 3-connected.
Suppose to its contrary, and let $\{u, v\}$ be a cutset of $G$.
Let $G_1$ and $G_2$ be two subgraphs of $G$ such that $V(G_1)\cap V(G_2)=\{u, v\}$ and $E(G_1)\cup E(G_2)=E(G)$, and let $f_i$ be a 3-coloring of $G_i$ for $i\in [2]$.  Then, $u\not\sim v$, and we may suppose by symmetry that $f_1(u)=f_1(v)=1$  and $f_2(u)=1$ and $f_2(v)=2$.

If $G_1$ has no $(1, i)$-$uv$-path for some $i\in\{2, 3\}$, we can modify $f_1$ to a 3-coloring $f'_1$ of $G_1$ with $f'_1(u)=1$ and $f'_1(v)=2$, and then $f'_1$ and $f_2$ form a 3-coloring of $G$, contradicting $\chi(G)=4$. For $i\in\{2, 3\}$, let $P_i$ be a shortest $(1, i)$-$uv$-path in $G_1$. Then, both $\l(P_2)$ and $\l(P_3)$ are of even. Since $G$ has no quadrilaterals, we may suppose that $\l(P_2)\ge 4$.
Let $P'$ be a shortest $(1, 2)$-$uv$-path in $G_2$ (such a path must exist). It is certain that $\l(P')$ is an odd number of at least 3 as $u\not\sim v$. Now, $P_2\cup P'$ is a big-odd-hole. This contradiction shows that $G$ is 3-connected.

\medskip

To prove $(b)$, let $S=\{u, u_1, u_2, u_3\}$ for some $u_1, u_2, u_3\in N(u)$, and assume that $S$ is a cutset of $G$. Let $H_1$ and $H_2$ be two subgraphs of $G$ such that $V(H_1)\cap V(H_2)=S$ and $E(H_1)\cup E(H_2)=E(G)$, and let $\phi_i$ be a 3-coloring of $H_i$ for $i\in [2]$ that maximize $|\{x\; :\; x\in S \mbox{ and } \phi_1(x)=\phi_2(x)\}|$. Without loss of generality, we may suppose that $\phi_1(u)=\phi_2(u)=1$, $\phi_1(u_1)=\phi_2(u_1)=2$,  $\phi_1(u_2)=2$, and $\phi_2(u_2)=3$. Before proving $(b)$, we need the following two claims on $\phi_1$ and $\phi_2$.

\begin{claim}\label{clm-con}
Let $x, y$ be two vertices of $S\setminus\{u\}$ with $\phi_1(x)=\phi_2(x)$ and $\phi_1(y)\neq \phi_2(y)$. Then, one of $H_1$ and $H_2$ has no $(2, 3)$-$xy$-paths.
\end{claim}
\pf Without loss of generality, we may suppose that $x=u_1$ and $y=u_2$.
If the claim is not true, let $P_i$, for $i\in [2]$,  be a shortest $(2, 3)$-$u_1u_2$-path in $H_i$. Now, $\l(P_1)\ge 4$ is even as $u_1, u_2\in N(u)$, and $\l(P_2)\ge 3$ is odd.

If $u_3\not\in V(P_1)\cup V(P_2)$, then $P_1\cup P_2$ is a big-odd-hole. If $u_3\in V(P_1)\cap V(P_2)$, let $P'_i=P_i[u_1, u_3]$ ($i\in[2]$) if $\phi_1(u_3)\neq \phi_2(u_3)$, and let $P'_i=P_i[u_2, u_3]$ ($i\in[2]$) if $\phi_1(u_3) = \phi_2(u_3)$, then  $P'_1\cup P'_2$ is a big-odd-hole.

So, we suppose, by symmetry, that $u_3\in V(P_1)\setminus V(P_2)$. If $u_3$ has no neighbors on $P_2$, then $P_1\cup P_2$ is a big-odd-hole. Otherwise, let $z_i$, for $i\in [2]$, be the neighbor of $u_3$ on $P_2$ which is closest to $u_i$. Now, either $P_1[u_1, u_3]\cup P_2[u_1, z_1]\cup z_1u_3$ (whenever $\phi_1(u_3)\neq \phi_2(u_3)$) or $P_1[u_2, u_3]\cup P_2[u_2, z_2]\cup z_2u_3$ (whenever $\phi_1(u_3) = \phi_2(u_3)$) is a big-odd-hole. This proves Claim~\ref{clm-con}. \qed

\begin{claim}\label{clm-con-1}
Let $x$ be vertex of $S$ with $\phi_1(x)\neq \phi_2(x)$, and let $\{y, z\}=S\setminus\{u, x\}$. Then, for each $i\in [2]$, $H_i$ has either a  $(2, 3)$-$xy$-path or a $(2, 3)$-$xz$-path.
\end{claim}
\pf Without loss of generality, we suppose that $x=u_2$, and suppose that $H_1$ has neither  $(2, 3)$-$u_1u_2$-paths nor  $(2, 3)$-$u_2u_3$-paths. Then, we can modify $\phi_1$ to a 3-coloring $\phi'_1$ of $H_1$ such that $\phi'_1(u_2)=3$ and $\phi'_1(w)= \phi_1(w)$ for other vertex $w$ of $S$. But now $\phi'_1$ and $\phi_2$ identify on more vertices of $S$ than those of $\phi_1$ and $\phi_2$, contradicting the choice of $\phi_1$ and $\phi_2$. \qed

\medskip

Now, we can prove ($b$). Without loss of generality,  we may suppose, by Claim~\ref{clm-con}, that $H_2$ has no $(2, 3)$-$u_1u_2$-path. Then, $H_2$ has $(2, 3)$-$u_2u_3$-paths by Claim~\ref{clm-con-1}, and consequently has no $(2, 3)$-$u_1u_3$-paths. Now, $H_1$ has no $(2, 3)$-$u_2u_3$-path by Claim~\ref{clm-con} which implies that $H_1$ has $(2, 3)$-$u_1u_2$-paths by Claim~\ref{clm-con-1}, and consequently has no $(2, 3)$-$u_1u_3$-path as $H_1$ has no $(2, 3)$-$u_2u_3$-paths. Note that $\{\phi_1(u_3), \phi_2(u_3)\}\subseteq \{2, 3\}$. We can modify $\phi_1$ to $\phi'_1$ such that $\phi'_1(u_3)\neq \phi_2(u_3)$ and $\phi'_1(x)=\phi_2(x)$ for $x\in S\setminus\{u_3\}$, and then modify $\phi_2$ to $\phi'_2$ that identifies with $\phi'_1$ on $S$, contradicting $\chi(G)=4$. This contradiction proves $(b)$.

\medskip

Let us now turn to prove $(c)$. Let $uv$ be an edge and $w$ be a vertex other than $u$ and $v$, and suppose to the contrary of $(c)$ that $\{u, v, w\}$ is a cutset. From $(b)$, we may suppose that $w\not\sim u$ and $w\not\sim v$. Let $G_1$ and $G_2$ be two subgraphs of $G$ with $V(G_1)\cap V(G_2)=\{u, v, w\}$ and $E(G_1)\cup E(G_2)=E(G)$. For $i\in[2]$, let $\phi_i$ be a 3-coloring of $G_i$. Since $\chi(G)=4$, we may suppose by symmetry that $|\{\phi_1(u), \phi_1(v), \phi_1(w)\}|=2$. Without loss of generality, suppose that
$$\phi_1(u)=\phi_2(u)=1, \mbox{ and } \phi_1(v)=\phi_2(v)=2.$$

We will show that
\begin{equation}\label{phi2-w}
\phi_2(w)=3.
\end{equation}

By symmetry, we suppose that $\phi_1(w)=1$, and suppose to the contrary of (\ref{phi2-w}) that $\phi_2(w)=2$.   Let $P_1$ be a shortest $(1, 3)$-$(u, w)$-path in $G_1$, $P'_1$ be a shortest $(1, 2)$-$(u, w)$-path in $G_1$, and $P_2$ be a shortest $(1, 2)$-$(u, w)$-path in $G_2$.
If $G_2$ has a $(2, 3)$-$(v, w)$-path, let $Q$ be such a path of shortest length. If $G_2$ has no $(2, 3)$-$(v, w)$-paths, let $\phi'_2$ be a 3-coloring  with $\phi'_2(u)=1, \phi'_2(v)=2$ and $\phi'_2(w)=3$, and let $Q'$ be a shortest $(1, 3)$-$(u, w)$-path, with respect to $\phi'_2$, in $G_2$ (such paths must exist since $\chi(G)=4$).

Since $P_1\cup P_2$ cannot be a big-odd-hole, we see  that
\begin{equation}\label{p1=2-v-in-p2}
\mbox{either $\l(P_1)=2$, or $v\in V(P_2)$ and $v$ has neighbors in $P_1-u$.}
\end{equation}

\begin{figure}[htbp]\label{fig-1}
  \begin{center}
    \includegraphics[width=6cm]{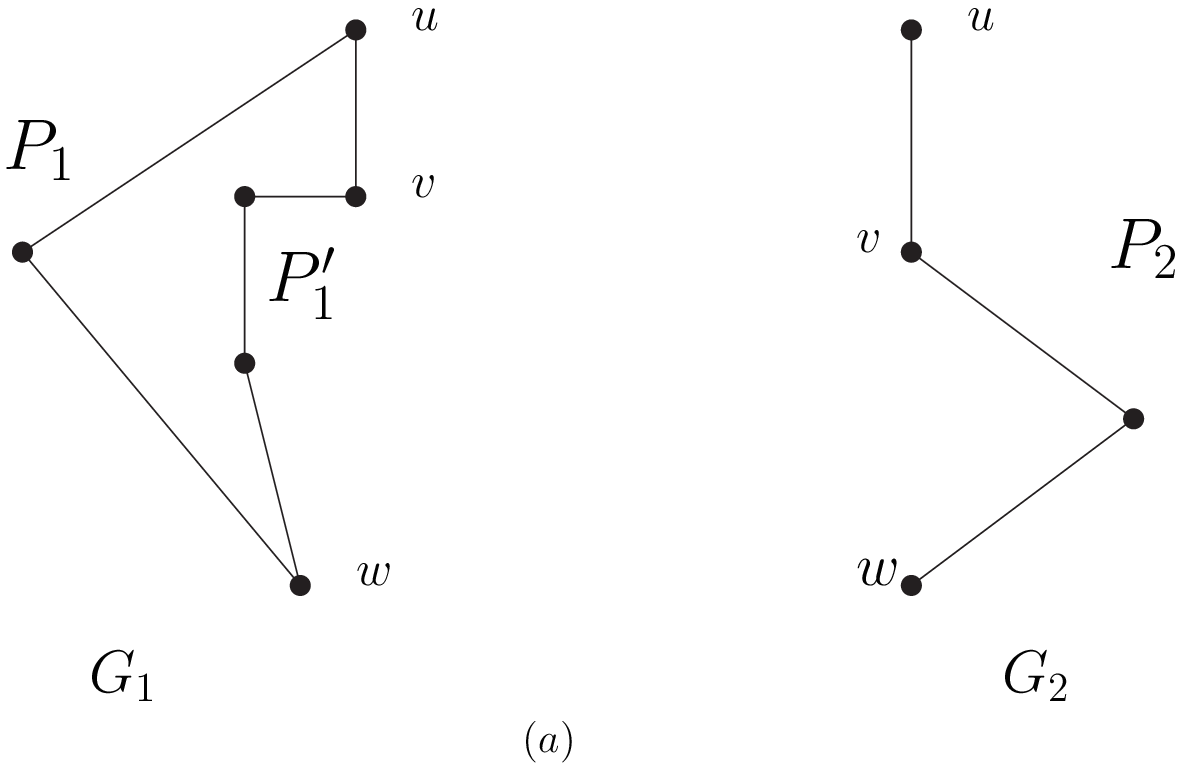} \hskip 40pt
    \includegraphics[width=6cm]{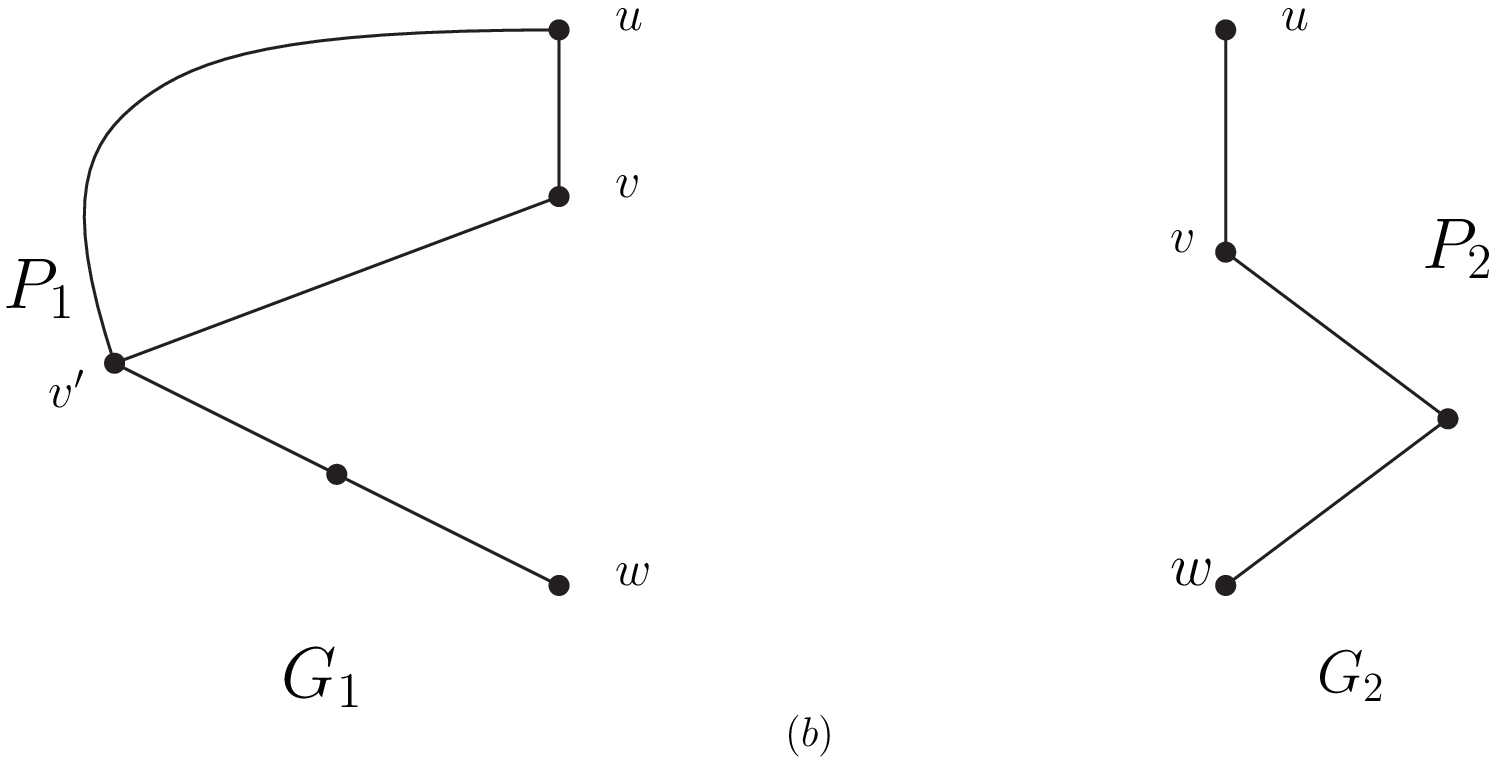}
  \end{center}
  \vskip -15pt
\caption{$\phi_1(u)=\phi_1(w)=\phi_2(u)=1$, $\phi_1(v)=\phi_2(v)=2$, and $\phi_2(w)=2$.}
\end{figure}

Suppose first that $\l(P_1)=2$. It is certain that $u$ is the unique neighbor of $v$ on $P_1$. If $\l(P_2)\ge 5$, then $P_1\cup P_2$ is a big-odd-hole. Therefore, $\l(P_2)=3$, and so $\l(P'_1)\ge 4$. If $v\not\in V(P'_1)$, then $v$ has no neighbors in $P'_1-u$, and so $P'_1\cup P_2$ is a big-odd-hole. If $v\in V(P'_1)$ and $\l(P'_1)\ge 6$, then $P'_1[v, w]\cup P_2$  induces a big-odd-hole.  If $v\in V(P'_1)$, $\l(P'_1)=4$, and $v\not\in V(P_2)$, then $P'_1\cup P_2$ is a big-odd-hole. Now, we have that $\l(P_1)=2$, $\l(P'_1)=4$ and $v\in V(P'_1)$,  and $\l(P_2)=3$ and $v\in V(P_2)$ (see Figure~\ref{fig-1}(a)).

If $G_2$ has a $(2, 3)$-$(v, w)$-path, then $\l(Q)\ge 4$ which implying that $P'_1[v, w]\cup Q$ is a big-odd-hole. So, $G_2$ has no $(2, 3)$-$(v, w)$-paths. If $\l(Q')=3$, then $v$ has no neighbor in $Q'-u$, and so $P'_1\cup Q'$ is a big-odd-hole. Otherwise, $\l(Q')\ge 5$, and so $P_1\cup Q'$ is a big-odd-hole. This contradiction shows that $\l(P_1)\ge 4$. Consequently, we have $v\in V(P_2)$ and $v$ has neighbors in $P_1-u$ by (\ref{p1=2-v-in-p2}).

\medskip

Let $v'$ be the neighbor of $v$ that is closest to $w$ on $P_1$. We claim first that
\begin{equation}\label{eqa-phiv'-1}
\phi_1(v')=1.
\end{equation}

If it is not the case, then $\phi_1(v')=3$.

Suppose that $v\in V(P'_1)$. We may suppose further that $P'_1=uvv_1v_2w$ and $P_2=uvw'w$ to avoid a big-odd-hole consisting of $P'_1-u$ and $P_2-u$. Consequently, $G_2$ has no $(2, 3)$-$(v, w)$-path as otherwise such a shortest path $Q$ together with $P'_1-u$ forms a big-odd-hole. Now, the path $Q'$ does exist and forms a big-odd-hole together with $P_1$. Therefore, $v\not\in V(P'_1)$, and so $u$ is the unique neighbor of $v$ in $P'_1$, which implies that $\l(P'_1)=2$ and $\l(P_2)=3$ to avoid a big-odd-hole consisting of $P'_1$ and $P_2$.
To avoid a $C_4$, we see that $\l(P_1[v', w])\ge 3$. To avoid the big-odd-hole $G[V(P_1)\cup V(Q')]$, we see that the path $Q$ does exist. If $u$ has no neighbor on $Q-v$, then $G[V(P_1)\cup V(Q)]$ is a big-odd-hole. Otherwise, let $u'$ be the neighbor of $u$ on $Q$ that is closest to $w$. Then, either $G[V(P_1)\cup V(Q[u', w])]$ or $G[V(P_1[v', w]\cup Q[u', w]\cup \{u, v\}]$ is a big-odd-hole. Therefore, (\ref{eqa-phiv'-1}) holds.

If $\l(P_1[v', w])\ge 5$, or $\l(P_2[v, w])\ge 4$, then $P_1[v', w]\cup P_2[v, w]$ is a big-odd-hole. We suppose so that $\l(P_1[v', w])=3$ and $\l(P_2[v, w])=2$, i.e, $G$ has a configuration as shown in Figure~\ref{fig-1}(b). With a similar argument as used to above Figure~\ref{fig-1}(a), one can deduce a contradiction by modifying $\phi_2$ to $\phi'_2$ with $\phi'_2(u)=1, \phi'_2(v)=2$ and $\phi'_2(w)=3$ and discussing a shortest $(1, 3)$-$(u, w)$-path with respect to $\phi'_2$. This proves (\ref{phi2-w}).

Next, we show that
\begin{equation}\label{phiW=must=1}
\mbox{if $\phi_1(u)=\phi_1(w)$ then $\psi(u)=\psi(w)$ in any 3-coloring $\psi$ of $G_1$.}
\end{equation}

Suppose that $\phi_1(w)=1$, and suppose that $G_1$ has another 3-coloring, say $\phi'_1$, such that $\phi'_1(u)=1$ and $\phi'_1(v)=\phi'_1(w)=2$. Since $\phi_2(w)=3$ by (\ref{phi2-w}), and since $\chi(G)=4$, we may choose $P_1$ to be a shortest $(1, 3)$-$(u, w)$-path in $G_1$ with respect to $\phi_1$, and choose $P'_1$ to be a shortest $(2, 3)$-$(v, w)$-path in $G_1$ with respect to $\phi'_1$. Let $P_2$ be a shortest $(1, 3)$-$(u, w)$-path in $G_2$, and let $P'_2$ be a shortest $(2, 3)$-$(v, w)$-path in $G_2$. Since $G$ has no big-odd-holes, we may suppose that $\l(P_1)=\l(P'_1)=2$ and $\l(P_2)=\l(P'_2)=3$.  Now, $P_2\cup P'_2+uv$ is a big-odd-hole, a contradiction. This proves (\ref{phiW=must=1}).

\medskip

Let $P_1$ be a shortest $(1, 2)$-$(u, w)$-path in $G_1$,  and let $P'_1$ be a shortest $(1, 3)$-$(u, w)$-path in $G_1$. Then, both $\l(P_1)$ and $\l(P'_1)$ are of even, and so one of them is at least 4 as $G$ has no $C_4$'s.

Since $\phi_2(w)=3$ by (\ref{phi2-w}), and since $\chi(G)=4$, we may choose $Q$ to be a shortest $(1, 3)$-$(u, w)$-path in $G_2$. Now, either $P'_1\cup Q$ is a big-odd-hole if $\l(P'_1)\ge 4$ or $\l(Q)\ge 5$, or $P_1\cup Q$ is a big-odd-hole if $\l(P'_1)=2$ and $\l(Q)=3$ (since $\l(P_1)\ge 4$ and $v$ has no neighbor on $Q$).  This proves $(c)$, and completes the proof of Theorem~\ref{4-critical}. \qed

\section{Proof of Theorem~\ref{coro-4-colorable-1}}

To prove Theorem~\ref{coro-4-colorable-1}, we in fact prove a little stronger conclusion that if a 3-connected graph $G\in {\cal G}_2$ has a 5-cycle $C$ and induces no 5-cycles sharing edges, then we can find an unstable 3-cutset of $G$ in $C$.

\begin{lemma}\label{no-two-5cycles}
Let $G\in {\cal G}_2$, and let $C=u_0u_1u_2u_3u_4u_0$ be a $5$-cycle in $G$. If $G$ is $3$-connected and induces no 5-cycles sharing edges, then it has a cutset of the form $\{u_i, u_{i+1}, w\}$ for some $i\in\{0, 1, 2, 3, 4\}$ and a vertex $w\in V(C)\setminus \{u_i, u_{i+1}\}$.
\end{lemma}
\pf Suppose to its contrary, and let $G$ be a counterexample. Then, $G$ is 3-connected, $\{u_i, u_{i+1}, w\}$ is not a cutset for each $i$ and each $w\in V(C)$.
	
Since $G$ has girth 5, and since $G$ has neither big-odd-hole nor 5-cycles sharing edges, we see that $L_1(C)$ is stable, no two vertices of $C$ have common neighbors in $L_1(C)$, and no two vertices of $L_1(C)$ have common neighbors in $L_2(C)$. Since $\delta(G)\ge 3$, we have $L_2(u_i)\cap L_2(C)\ne\emptyset$ for each $i\in\{0, 1, 2, 3, 4\}$.

Let $m$ be an integer in $\{0, 1, 2, 3, 4\}$.  Since $\{u_{m+1}, u_{m+2}, u_{m+4}\}$ is not a cutset of $G$, we see that $G-V(C)$ has an induced $(L_2(C)\cap L_2(u_m), L_2(C)\cap L_2(u_{m+3}))$-path. We will show that each such path has no neighbors of either $u_{m+1}$ or $u_{m+2}$.

\begin{claim}\label{eqa-Ui-Ui+1-00}
Let $Q$ be an induced $(L_2(C)\cap L_2(u_m), L_2(C)\cap L_2(u_{m+3})$-path in $G-V(C)$. Then,  $N_{Q^*}(u_{m+1})\cup N_{Q^*}(u_{m+2})=\emptyset$ and $N_{Q^*}(u_{m+4})\neq\emptyset$.
\end{claim}
\pf By symmetry, it suffices to prove the case that $m=0$. For $i\in \{1, 3\}$, let $x_i$ be the vertices of $Q$ in $L_2(u_i)$, and let $z_i$  be the common neighbor of $u_i$ and $x_i$. It is certain that $N_{Q^*}(u_0)\cup N_{Q^*}(u_{3})=\emptyset$, and both $z_0$ and $z_3$ have no neighbors in $Q^*$.

If $Q^*\cap L_1(C)=\emptyset$, then either $Qz_0u_0u_1u_2u_3z_3$ or $Qz_0u_0u_4u_3z_3$ is a big-odd-hole. So, we suppose that $Q^*\cap L_1(C)\ne \emptyset$, and let $v_0, v_1, \dots, v_r$ be the vertices in $Q^*\cap L_1(C)$ following the direction from $x_0$ to $x_{3}$.

To prove our claim, we need only to verify, by symmetry, that
\begin{equation}\label{eqa-U1-Q0}
N_{Q^*}(u_1)=\emptyset.
\end{equation}

Suppose to the contrary of  (\ref{eqa-U1-Q0}), and suppose,  without loss of generality, that $v_0\sim u_1$. It is certain that $v_0\not\sim u_4$, and $N_{Q^*}(u_3)=\emptyset$ by the definition of $Q$. We will deduce a contradiction that $Qz_3u_3u_4u_0z_0$ is a big-odd-hole.

Suppose that $N_{Q^*}(u_4) \neq \emptyset$. Let $j'\in [r]$ be the minimum integer such that $u_4\sim v_{j'}$, and let $j''\in\{0, 1, \ldots, j'-1\}$ be the largest such that $u_1\sim v_{j''}$. If $u_2$ has neighbors in $\{v_{j''+1}, \ldots, v_{j'-1}\}$, let $j_1\in \{j''+1, \ldots, j'-1\}$ be the largest such that $u_2\sim v_{j_1}$. Now, either $Q[v_{j''}, v_{j'}]u_4u_0u_1$ or $Q[v_{j''}, v_{j'}]u_4u_3u_2u_1$ is a big-odd-hole if $u_2$ has no neighbors in $\{v_{j''+1}, \ldots, v_{j'-1}\}$, and either $Q[v_{j_1}, v_{j'}]u_4u_0u_1u_2$ or $Q[v_{j_1}, v_{j'}]u_4u_3u_2$ is a big-odd-hole otherwise. Therefore,
\begin{equation}\label{eqa-Ui-Ui+1}
N_{Q^*}(u_4) = \emptyset.
\end{equation}

Since $G$ has neither big-odd-holes nor 5-cycles sharing edges, we have that, for each $i\ge 0$,
\begin{equation}\label{eqa-same}
\mbox{$Q[v_i, v_{i+1}]$ is odd if neither $u_1\in N(v_i)\cap N(v_{i+1})$ nor $u_2\in N(v_i)\cap N(v_{i+1})$}
\end{equation}
as otherwise $Q[v_i, v_{i+1}]$ together with $u_1$ and $u_2$ forms a big-odd-hole, and
\begin{equation}\label{eqa-distinct}
\mbox{$Q[v_i, v_{i+1}]$ is even or odd length 3 if either $u_1\in N(v_i)\cap N(v_{i+1})$ or $u_2\in N(v_i)\cap N(v_{i+1})$.}
\end{equation}

Now we let $j_0\in\{0, 1, \ldots, r\}$ be the largest such that $u_1\sim v_{j_0}$, and let $i_0\in\{0, 1, \ldots, r\}$  be the smallest such that $u_2\sim v_{i_0}$.   It is certain that $i_0>0$ as $u_1\sim v_0$, and $j_0<r$ as otherwise $Q[v_r, x_3]z_3u_3u_4u_0u_1$ or $Q[v_r, x_3]z_3u_3u_2u_1$ is a big-odd-hole. Since $N_{Q^*}(u_0)=N_{Q^*}(u_3)=\emptyset$, we have, by (\ref{eqa-Ui-Ui+1}) and (\ref{eqa-same}), that
\begin{equation}\label{eqa-U1-I0-1}
\mbox{$u_1\sim v_{i_0-1}$, $u_2\sim v_{j_0+1}$, and both $Q[v_{i_0-1}, v_{i_0}]$ and $Q[v_{j_0}, v_{j_0+1}]$ are odd paths}.
\end{equation}

We claim further that
\begin{equation}\label{eqa-I0-J0}
\mbox{$Q[v_{i_0}, v_{j_0}]$ must be odd.}
\end{equation}

Suppose to its contrary that $Q[v_{i_0}, v_{j_0}]$ is even. We have that $i_0< j_0-1$, as otherwise either $i_0=j_0+1$ or $i_0=j_0-1$, and in either case $Q[v_{i_0}, v_{j_0}]$ together with $u_1$ and $u_2$ forms a big-odd-hole.
Let $S=\{v_{i_0}, v_{i_0+1}, \ldots, v_{j_0}\}$. We order all the vertices of $S\cap N(u_1)$ into a sequence $S'=v_{i^{(1)}}v_{i^{(2)}}\ldots v_{i^{(t')}}$ along $Q$ from $x_0$ to $x_3$, and order all the vertices of $S\cap N(u_2)$ into a sequence $S''=v_{j^{(1)}}v_{j^{(2)}}\ldots v_{j^{(s')}}$ similarly. Since $u_1\sim v_{i_0-1}$, we see that  $Q[v_{i_0-1},v_{i^{(1)}}]$ must be  even. Since $Q[v_{i_0-1}, v_{i_0}]$ is odd  by (\ref{eqa-same}), we have that
\begin{equation}\label{eqa-odd-length-0}
\mbox{$Q[v_{i_0}, v_{i^{(1)}}]$ is odd,}
\end{equation}
and so
\begin{equation}\label{eqa-odd-length-1}
\mbox{$Q[v_{i^{(1)}}, v_{j_0}]$ is odd}
\end{equation}
as $Q[v_{i_0}, v_{j_0}]$ is even by our assumption.

Now, let us check the number of 5-cycles in $G[V(Q[v_{i_0}, v_{j_0}])\cup\{u_i\}]$ for each $i\in [2]$. We will deduce a contradiction by counting the 5-cycles in two ways, and thus prove  (\ref{eqa-I0-J0}).

For each integer $h\in\{1,\dots,t'-1\}$, we have three possibilities about the segment $Q[v_{i^{(h)}},v_{i^{(h+1)}}]$: $(A)$ $i^{(h+1)}=i^{(h)} + 1$ and $Q[v_{i^{(h)}}, v_{i^{(h+1)}}]$ is odd, i.e., $Q[v_{i^{(h)}}, v_{i^{(h+1)}}]u_1$ is a 5-cycle; $(B)$ $i^{(h+1)}=i^{(h)} + 1$ and $Q[v_{i^{(h)}}, v_{i^{(h+1)}}]$ is even; and $(C)$ $i^{(h+1)} > i^{(h)} + 1$, and  $Q[v_{i^{(h)}}, v_{i^{(h+1)}}]$ is even.

It is easy to check that, $Q[v_{i^{(h)}}, v_{i^{(h+1)}}]u_1$ is a 5-cycle if and only if $(A)$ occurs. Since $Q[v_{i^{(1)}}, v_{j_0}]$ is odd by (\ref{eqa-odd-length-1}), we see that the number of segments satisfying $(A)$ must be odd, which implies that there are  odd number of 5-cycles in $G[V(Q[v_{i_0}, v_{j_0}])\cup\{u_1\}]$ by (\ref{eqa-distinct}).  So, we have, by symmetry, that
\begin{equation}\label{eqa-odd-5cycles}
\mbox{there are odd number of 5-cycles in $G[V(Q[v_{i_0}, v_{j_0}])\cup\{u_i\}]$ for each $i\in [2]$}.
\end{equation}

For $i\in [2]$, we call the 5-cycles in  $G[V(Q[v_{i_0}, v_{j_0}])\cup\{u_i\}]$ as $Ui$-5-cycles. It is clear that no 5-cycle can be  both $U1$-5-cycle and $U2$-5-cycle.  Next, let us count the number of $U2$-5-cycles in $G[V(Q[v_{i_0}, v_{j_0}])\cup\{u_2\}]$.

Suppose the third case $(C)$ occurs for some $h$. Then, all vertices of $Q[v_{i^{(h)}},v_{i^{(h+1)}}]\cap S$ are neighbors of $u_2$. By (\ref{eqa-distinct}), both $Q[v_{i^{(h)}}, v_{i^{(h)}+1}]$ and $Q[v_{i^{(h+1)}-1}, v_{i^{(h+1)}}]$ are odd, and so the number of $U2$-5-cycles in $G[V(Q[v_{i^{(h)}}, v_{i^{(h+1)}}])\cup\{u_2\}]$ must be even as otherwise $Q[v_{i^{(h)}}, v_{i^{(h+1)}}]u_1$ is a big-odd-hole.

If $i^{(1)}-1=i_0$, then we can find totally even number of $U2$-5-cycles in $G[V(Q[v_{i_0}, v_{j_0}])\cup\{u_2\}]$. If $i^{(1)}-1>{i_0}$, then $Q[v_{i_0}, v_{i^{(1)}-1}]$ is even because both $Q[v_{i^{(1)}-1}, v_{i^{(1)}}]$ and $Q[v_{i_0}, v_{i^{(1)}}]$ are odd by (\ref{eqa-same}) and (\ref{eqa-odd-length-0}), and so $G[V(Q[v_{i_0}, v_{i^{(1)}-1}])\cup\{u_2\}]$ must have even number of $U2$-5-cycles. In either cases, we have  even number of $U2$-5-cycles in $G[V(Q[v_{i_0}, v_{j_0}])\cup\{u_2\}]$, contradicting (\ref{eqa-odd-5cycles}). This proves (\ref{eqa-I0-J0}).

\medskip

Note that $\l(Q)=\l(Q[x_0, v_{i_0}])+\l(Q[v_{i_0}, v_{j_0}])+\l(Q[v_{j_0}, x_3])$. Since both $u_0$ and $u_3$ have no neighbor in $Q^*$, we see that $Q[v_{j_0},x_3]$ must be even as otherwise $Q[v_{j_0}, x_3]z_3u_3u_4u_0u_1$ is a big-odd-hole. Similarly, $Q[x_0, v_{i_0}]$ is even.  But then, $Q$ is odd as $Q[v_{i_0}, v_{j_0}]$ is odd by (\ref{eqa-I0-J0}), and so $Qz_3u_3u_4u_0z_0$ is a big-odd-hole. This contradiction proves (\ref{eqa-U1-Q0}), and thus proves Claim~\ref{eqa-Ui-Ui+1-00}. \qed

\medskip

By  taking $m=2$ in Claim~\ref{eqa-Ui-Ui+1-00}, we see that $G-V(C)$ has an induced $(L_2(C)\cap L_2(u_0), L_2(C)\cap L_2(u_2))$-path, say $P$, such that no vertex of $V(C)\setminus \{u_1\}$ has neighbors in $P^*$. To forbid a big-odd-hole, $V(P)\cap L_1(C)\neq \emptyset$.
Let $w$ be a vertex of $P^*\cap N(u_1)$, and let $y_1$ be a neighbor of $w$ in $P$. It is certain that $y_1\in L_2(C)\cap L_2(u_1)$.

By our assumption, $G$ has no cutset of the form $\{u_m, u_{m+1}, w\}$ for any $0\le m\le 4$ and any $w\in V(C)$. We see that $G-\{u_0, u_1, u_2\}$ is connected, and thus contains a path from $y_1$ to $L_2(C)\cap (L_2(u_3)\cup L_2(u_4))$. Let $P'$ be such a path of shortest length. It is certain that $P'$ is contained in $G-V(C)$ and $(P')^*\cap (L_1(u_3)\cup L_1(u_4))=\emptyset$. Without loss of generality, we suppose that $P'$ is a $(y_1, L_2(C)\cap L_2(u_3))$-path. Now, we can find, in $P\cup P'$,  an induced $(L_2(C)\cap L_2(u_0), L_2(C)\cap L_2(u_3))$-path, say $Q$, with $N_{Q^*}(u_4)=\emptyset$. This contradiction to Claim~\ref{eqa-Ui-Ui+1-00} proves Lemma~\ref{no-two-5cycles}. \qed

\medskip

It is certain that Theorem~\ref{coro-4-colorable-1} follows directly from Lemma~\ref{no-two-5cycles}.


\section{Proof of Theorem~\ref{4-critical1}}

This section is devoted to the proof Theorem~\ref{4-critical1}: Let $G$ be a $3$-connected graph in ${\cal G}_2$ without unstable 3-cutset. Suppose that $G$ is not isomorphic to the Petersen graph. Then, $G$ does not induce $\theta^+$.

We need some new notations.
Let $H$ be a proper induced subgraph of $G$. An {\em ear} $F$ of $H$ in $G$ is a path connecting two vertices of $H$ with $F^*\subseteq V(G)\setminus V(H)$. Let $F$ be an ear of $H$. The two common vertices of $F$ and $H$ are called the {\em attachments} of $F$ in $H$. If $F$ itself is an induced path, then we call $F$ an {\em induced ear}. Recall that for two induced subgraph $F_1$ and $F_2$ of $G$, we simply write $G[V(F_1)\cup V(F_2)]$ as $F_1\cup F_2$.

Let $F$ be an induced ear of $H$ with attachments $\{x, y\}$. We say that $F$ is a {\em strong induced ear} of $H$ if $E(H\cup F)\setminus (E(H)\cup E(F))=E(N_H(x)\cap N_H(y), F^*)$, i.e., the only possible edges of $H\cup F-E(H)\cup E(F)$ are those between the common neighbor of $x$ and $y$ in $H$ and the internal vertices of $F$.

Below Proposition~\ref{ear} is very useful in our proof.

\begin{proposition}\label{ear}
Let $G$ be a $3$-connected graph in ${\cal G}_2$, let $H$ be a proper induced subgraph of order at least $3$ such that each vertex of $V(G)\setminus V(H)$ has at most one neighbor in $H$. Then, $H$ has a strong induced ear of length at least $3$ with nonadjacent attachments.
\end{proposition}
\pf Since $G$ is 3-connected and has no $C_3$, we see that every vertex of $V(G)\setminus V(H)$ has three internally disjoint paths to $H$, and we may choose $P$ to be a shortest path connecting two nonadjacent vertices of $H$ with $P^*\subseteq V(G)\setminus V(H)$. Let $x_1$ and $x_2$ be the ends of $P$.
Since each vertex of $V(G)\setminus V(H)$ has at most one neighbor in $H$, it is clear that $\l(P)\ge 3$. Let $S=\{x_1, x_2\}\cup (N(x_1)\cap N(x_2))$. The minimality of $\l(P)$ shows that no vertex of $V(H)\setminus S$ has neighbor in $P^*$. This proves Proposition~\ref{ear}. \qed

We use ${\cal P}$ to denote the Petersen graph, use ${\cal P}^-$ to denote the graph obtained from ${\cal P}$ by removing a vertex. Recall that $\theta^+$ is the graph obtained from ${\cal P}$  by removing two adjacent vertices (see Figure~\ref{fig-3}). From now on, let $G$ be a $3$-connected graph in ${\cal G}_2$ of which each 3-cutset is stable.

\medskip

Before proving that $G$ induces no $\theta^+$, we show firstly that $G$ induces no ${\cal P}^-$.

\begin{proposition}\label{petersen}
$G$ induces no ${\cal P}^-$.
\end{proposition}
\pf First we show that $G$ does not induce ${\cal P}$. Suppose to its contrary, let $H$ be an induced ${\cal P}$ of $G$ with $V(H)=\{u_i, v_i : i=0,1,2,3,4\}$ and $E(H)=\{u_iu_{i+1}, v_iv_{i+2}, u_iv_i: i=0,1,2,3,4\}$.
It is certain that $G\neq H$, and  every vertex in $V(G)\setminus V(H)$ has at most one neighbor in $H$ since any two vertices of $H$ are contained in a 5-cycle.  By Proposition~\ref{ear}, we may suppose that $H$ has a strong induced ear of length at least 3, say $P_0$, with attachments $\{u_0, u_2\}$. Then, $e(H-\{u_0, u_1, u_2\}, P_0-\{u_0, u_2\})=0$, and so either $P_0u_0u_4u_3u_2$ or $P_0u_0u_4v_4v_2u_2$ is a big-odd-hole. Therefore, $G$ does not induce ${\cal P}$.

Next, we show that $G$ does not induce the Petersen graph minus an edge. Suppose to its contrary, let $H'$ be an induced subgraph of $G$ with  $V(H')=\{u_i, v_i :  i=0,1,2,3,4\}$ and $E(H')=\{u_iu_{i+1}, v_iv_{i+2}, u_iv_i: i=0,1,2,3,4\}\setminus\{u_0v_0\}$. If $V(G)\setminus V(H')$ has a vertex, say $x$, with two neighbors in $H'$, then $N_{H'}(x)\in \{\{u_0, v_0\}, \{u_0, v_3\}, \{u_1, v_0\}, \{u_4, v_0\}\}$, and one can easily find a big-odd-hole in $G$. Therefore, $|N_{H'}(x)|\le 1$ for every vertex $x\not\in H'$. By Proposition~\ref{ear}, $H'$ has a strong induced ear of length at least 3, say $P_1$, with nonadjacent attachments, say $\{x, y\}$.  Note that in ${\cal P}$, every pair of nonadjacent vertices is connected by two edge disjoint induced paths of even length at least 4, and connected by two edge disjoint induced paths of odd length. In $H'$, every pair of nonadjacent vertices  is connected by an induced path of even length at least 4, and by an induced path of odd length. Again, a big-odd-hole occurs in $H'\cup P_1$. This contradiction shows that $G$ does not induce the Petersen graph minus an edge.

\medskip

Finally, suppose to the contrary of Proposition~\ref{petersen} that $G$ induces a ${\cal P}^-$, say $H''$. One can easily check that in $H''$, every pair of nonadjacent vertices  is connected by an induced path of even length at least 4, and connected by an induced path of odd length.  Since $G$ induces neither ${\cal P}$ nor ${\cal P}^-$, each vertex of $V(G)\setminus V(H'')$ has at most one neighbor in $H''$. By Proposition~\ref{ear}, $H''$ has a strong induced ear of length at least 3  with nonadjacent attachments. Then, a big-odd-hole occurs in $G$ (we leave the details to readers), a contradiction.  \qed

\bigskip

Now let us turn to the proof of $G$ inducing no $\theta^+$.

Suppose to its contrary, let $H$ be an induced  $\theta^+$ in $G$ consisting of a cycle $v_0v_1v_2v_3v_4v_5v_6v_7$ together with chords $v_0v_4$ and $v_2v_6$ (see Figure~\ref{fig-3}). By Proposition~\ref{petersen}, one can verify easily that
\begin{equation}\label{theta+1}
|N_{H}(w)|\le 1 \mbox{ for every vertex $w\in V(G)\setminus V(H)$}.
\end{equation}

By Proposition~\ref{ear} and (\ref{theta+1}), $H$  must have strong induced ears of length at least 3 with nonadjacent attachments. We call a strong induced ear of $H$ as a $(v_i, v_{i+2})$-{\em bad ear} if it has length 3 and its attachments are $v_i$ and $v_{i+2}$ for some $i\in\{1, 3, 5, 7\}$.

\begin{claim}\label{bad-ear-0}
Every strong induced ear of $H$ with nonadjacent attachments is a bad ear.
\end{claim}
\pf Let $P$ be a strong induced ear of $H$ with nonadjacent attachments, say $\{x, y\}$. It is certain that $\l(P)\ge 3$ by (\ref{theta+1}). To prove our claim, we need to verify that
$\{x, y\}=\{v_i, v_{i+2}\}$ for some $i\in\{1, 3, 5, 7\}$, and $\l(P)=3$. Here the summation of subindexes are taken modulo 7.

Suppose that $\{x, y\}\ne \{v_i, v_{i+2}\}$ for any $i\in\{1, 3, 5, 7\}$. Then, we may assume by symmetry that $\{x, y\}\in\{\{v_0, v_2\}, \{v_0, v_3\}, \{v_1, v_5\}\}$.

If $\{x, y\}=\{v_0, v_2\}$ let $C=v_0v_4v_3v_2P$ (whenever $\l(P)$ is even) or $C=v_0v_4v_5v_6v_2P$ (whenever $\l(P)$ is odd), if $\{x, y\}=\{v_0, v_3\}$ let $C=v_0v_1v_2v_3P$ (whenever $\l(P)$ is even) or $C=v_0v_7v_6v_2v_3P$ (whenever $\l(P)$ is odd), if $\{x, y\}=\{v_1, v_5\}$  let $C=v_1v_2v_6v_5P$ (whenever $\l(P)$ is even) or $C=v_1v_2v_3v_4v_5P$ (whenever $\l(P)$ is odd), then $C$ is a big-odd-hole in each case.

Therefore, $\{x, y\}=\{v_i, v_{i+2}\}$ for some $i\in\{1, 3, 5, 7\}$. Suppose  by symmetry that $\{x, y\}=\{v_1, v_3\}$, and suppose that $\l(P)>3$. Then, $N_{P^*}(v_2)\neq \emptyset$
as otherwise either $Pv_1v_2v_3$ or $Pv_1v_0v_4v_3$ is a big-odd-hole. By our assumption, $G-\{v_1, v_2, v_3\}$ is connected. Let $P'$ be a shortest $(\{v_0, v_4, v_5, v_6, v_7\}, P^*)$-path in $G-\{v_1, v_2, v_3\}$, let $x$ be the end of $P'$ on $P$, and let $y$ be the other end of $P'$. We may suppose that $y\in \{v_4, v_5, v_6\}$ by symmetry. If $y=v_4$, let $P''$ be a shortest $v_2y$-path in $P^*\cup P'\cup \{v_2\}$, then either $v_2v_6v_5v_4P''$ or $v_2v_6v_7v_0v_4P''$ is a big-odd-hole. If $y\in \{v_5, v_6\}$, let $P''$ be a shortest $v_1y$-path in $P[v_1, x]\cup P'$, then $v_1v_0v_4v_5P''$ or $v_1v_0v_7v_6v_5P''$ is a big-odd-hole (whenever $y=v_5$), and $v_1v_0v_7v_6P''$ or $v_1v_0v_4v_5v_6P''$ is a big-odd-hole (whenever $y=v_6$). This contradiction proves Claim~\ref{bad-ear-0}.\qed

\medskip

From Proposition~\ref{ear} and Claim~\ref{bad-ear-0}, $H$ must have bad ears. Since $|N_H(t)|\le 1$ for each vertex $t\not\in H$ by (\ref{theta+1}), and since $G$ induces no big-odd-hole,  we see that at most two pairs from $\{\{v_1, v_3\}, \{v_1, v_7\}, \{v_3, v_5\}, \{v_5, v_7\}\}$ may be the attachments of bad ears. We suppose, by symmetry, that
\begin{equation*}
\mbox{$H$ has $(v_1, v_3)$-bad ears, but no  $(v_3, v_5)$-bad ears.}
\end{equation*}

Next, we study the subgraphs consisting of $H$ and all its bad ears. Let $H'$ be a maximal induced subgraph of $G$ such that $V(H)\subseteq V(H')$ and each vertex of $V(H')\setminus V(H)$ is in a bad ear of $H$.

\begin{claim}\label{bad-ears}
Each vertex of $V(H')\setminus V(H)$ $($resp. each edge of $E(H')\setminus E(H))$ is in a unique bad ear.
\end{claim}
\pf If $H$ has two bad ears, say $P_1$ and $P_2$, sharing a common internal vertex, say $w$, we may suppose by symmetry that $P_1=v_1wxv_7$ and $P_2=v_1wyv_3$ (as $|N_H(t)|\le 1$ for each vertex $t\not\in H$ by (\ref{theta+1})), then $v_2v_3ywxv_7v_6v_2$ is a big-odd-hole. Therefore, no two distinct bad ears share internal vertex.

Since $G$ has girth 5, there is no edges between two bad ears with the same attachments.

To prove the second statement, we assume to its contrary that $xy$ is an edge between two distinct bad ears, say $P_1$ and $P_2$, with $x\in V(P_1)$ and $y\in V(P_2)$. Without loss of generality, suppose that $P_1=v_1ww'v_3$. Then, $P_2$ is either a $v_1v_7$-path, or a $v_5v_7$-path. If $P_2=v_1zz'v_7$, then $x=w'$ and $y=z'$ which implies that   $w'v_3v_4v_5v_6v_7z'w'$ is a big-odd-hole. If $P_2=v_5zz'v_7$, we suppose by symmetry that $x=w$, then either $ww'v_3v_4v_5zz'w$ (whenever  $y=z'$) or $ww'v_3v_2v_6v_5zw$ (whenever $y=z$) is a big-odd-hole. This completes the proof of Claim~\ref{bad-ears}. \qed

\begin{figure}[htbp]\label{fig-4}
  \begin{center}
    \includegraphics[width=5cm]{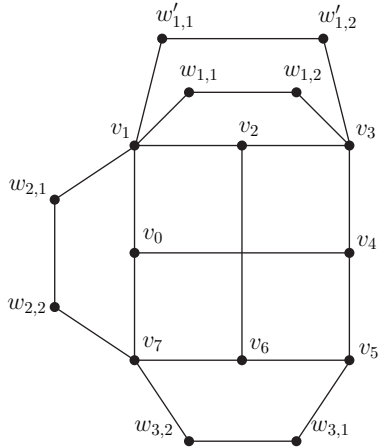}
  \end{center}
  \vskip -15pt
\caption{$H$ and some possible bad ears.}
\end{figure}

\begin{claim}\label{bad-ears-1}
Each vertex in $V(G)\setminus V(H')$ has at most one neighbor in $H'$.
\end{claim}
\pf Suppose that the claim does not hold, and let $x$ be a vertex in $V(G)\setminus V(H')$ with two neighbors, say $z$ and $z'$, in $H'$. By (\ref{theta+1}), it is certain that $|\{z, z'\}\cap V(H)|\le 1$. We use $P_1=v_1w_{1, 1}w_{1,2}v_3$, $P'_1=v_1w'_{1, 1}w'_{1,2}v_3$,  $P_2=v_1w_{2, 1}w_{2,2}v_7$, and $P_3=v_5w_{3, 1}w_{3,2}v_7$ to represent four possibly existing bad ears. See Figure~\ref{fig-4}.

Suppose that $|\{z, z'\}\cap V(H)|= 1$, say $z\in V(H)$ by symmetry. We may suppose further that $z\in\{v_1, v_2\}$. If $z=v_1$, then $z'\in \{w_{3, 1}, w_{3, 2}\}$, and either $v_1v_2v_3v_4v_5w_{3, 1}xv_1$ (when $z'=w_{3, 1}$) or $v_1v_2v_6v_5w_{3, 1}w_{3, 2}xv_1$ (when $z'=w_{3, 2}$) is a big-odd-hole. If $z=v_2$, we may suppose by symmetry that $z'\in \{w_{2, 2}, w_{3, 2}\}$, and then either $xv_2v_3v_4v_0v_7w_{2, 2}x$ (when $z'=w_{2, 2}$) or $xv_2v_3v_4v_0v_7w_{3, 2}x$ (when $z'=w_{3, 2}$) is a big-odd-hole.

Therefore, $|\{z, z'\}\cap V(H)|= 0$. By symmetry, we suppose that $z\in\{w_{1, 1}, w_{1, 2}\}$.

If $z$ and $z'$ are on two distinct ears with the same attachments, we may suppose by symmetry that $z=w_{1, 1}$ and $z'=w'_{1, 2}$ (as $G$ has girth 5), then  $v_1w_{1, 1}xw'_{1, 2}v_3v_4v_0v_1$ is a big-odd-hole.
So, $z$ and $z'$ cannot be on the ears with the same attachments.

If $z=w_{1, 1}$, then $z'\in\{w_{2, 2}, w_{3, 1}, w_{3, 2}\}$, and a big-odd-hole $C'$ appears, where $C'=xw_{1, 1}v_1v_2v_6v_7w_{2, 2}x$ whenever $z'=w_{2, 2}$, $C'=xw_{1, 1}v_1v_2v_6v_5w_{3, 1}x$  whenever $z'=w_{3, 1}$, and  $C'=xw_{1, 1}v_1v_2v_6v_7w_{3, 2}x$  whenever $z'=w_{3, 2}$. One can deduce the same contradiction if $z=w_{1, 2}$. These contradictions prove Claim~\ref{bad-ears-1}. \qed

\medskip

Recall that $H$ has at most two families of bad ears of which each family has the same attachments, and we always suppose that $H$ has $(v_1, v_3)$-bad ears and has no $(v_3, v_5)$-bad ears. Let $v_1w_1w_2v_3$ be a bad ear. We distinguish  two cases upon whether $H$ has $(v_5, v_7)$-bad ears or not.

\medskip

\noindent{\bf Case} 1. We suppose first that $H$ has no $(v_5, v_7)$-bad ears.

Recall that every 3-cutset of $G$ must be a stable set. Let $P$ be a shortest $(v_5, V(H)\setminus\{v_4, v_5, v_6\})$-path in $G-\{v_4, v_6\}$, and let $x$ be the end of $P$ other than $v_5$. By symmetry, we may suppose that $x\in \{v_0, v_1, v_7\}$. Then, $\l(P)\ge 3$ by (\ref{theta+1}). Since $H$ has no strong induced ears with attachments $\{v_4, v_6\}$ by Claim~\ref{bad-ear-0}, we see that
\begin{equation}\label{v4-v6}
\mbox{either $v_4$ or $v_6$ has no neighbors on $P-v_5$}.
\end{equation}

\begin{figure}[htbp]\label{fig-5}
  \begin{center}
    \includegraphics[width=10cm]{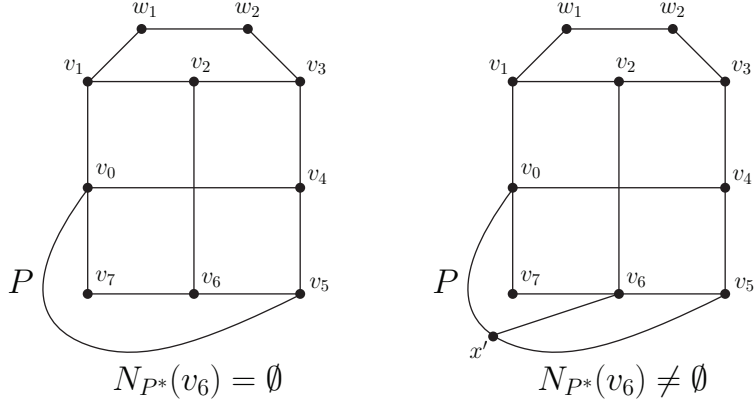}
  \end{center}
  \vskip -15pt
\caption{$H$ has no $(v_5, v_7)$-bad ears and  $P$ is a $v_0v_5$-path.}
\end{figure}

First, suppose that $x=v_0$ (see Figure~\ref{fig-5}). If $v_6$ has no neighbors in $P-v_5$, then either $v_5v_6v_7v_0P$ or $v_5v_6v_2v_1v_0P$  is a big-odd-hole. Otherwise, let $x'$ be the neighbor of $v_6$ on $P-v_5$ which is  closest to $x$, and let $P'=v_6x'P[x', x]$. Now, $v_4$ has no neighbors on $P'$ by (\ref{v4-v6}), and thus either $v_6v_2v_1v_0P'$  or $v_6v_2v_3v_4v_0P'$  is a big-odd-hole, a contradiction.

Next, suppose that $x=v_1$. If $v_4$ has no neighbors in $P-v_5$, then either $v_5v_4v_0v_1P$ or $v_5v_4v_3v_2v_1P$  is a big-odd-hole. Otherwise, let $x'$ be the neighbor of $v_4$ on $P-v_5$ which is closest to $x$, and let $P'=v_4x'P[x', x]$. Now, $v_6$ has no neighbors on $P'$ by (\ref{v4-v6}), and thus either $v_4v_3v_2v_1P'$  or $v_4v_5v_6v_2v_1P'$ is a big-odd-hole, a contradiction.

Finally, we suppose that $x=v_7$. By our assumption that $H$ has no $(v_5, v_7)$-bad ears, $v_4$ must have neighbors in $P-v_5$. Let $x'$ be the neighbor of $v_4$ on $P-v_5$ which is closest to $x$, and let $P'=v_4x'P[x', x]$. Then, $v_6$ has no neighbors on $P'$ by (\ref{v4-v6}), and thus either $v_4v_5v_6v_7P'$  or $v_4v_3v_2v_6v_7P'$  is a big-odd-hole, a contradiction again.

\medskip

\begin{figure}[htbp]\label{fig-6}
  \begin{center}
    \includegraphics[width=4cm]{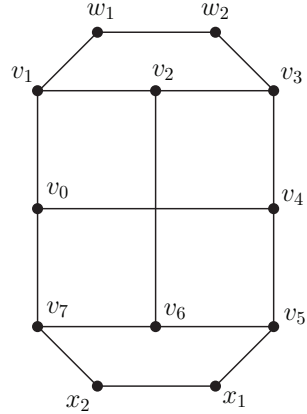}
  \end{center}
  \vskip -15pt
\caption{$H$ together with a $(v_1, v_3)$-bad ear and a $(v_5, v_7)$-bad ear}
\end{figure}

\noindent{\bf Case} 2. Now, we suppose that $H$ has $(v_5, v_7)$-bad ears, and let $v_5x_1x_2v_7$ be a bad ear (see Figure~\ref{fig-6}). Then, $H$ has no bad ears with attachments $\{v_1, v_7\}$ or $\{v_3, v_5\}$.

First, we claim that any four consecutive vertices on the cycle $v_0v_1v_2v_3v_4v_5v_6v_7$ cannot be a cutset. To be precisely, let ${\cal S}=\{\{v_i, v_{i+1}, v_{i+2}, v_{i+3}\}\;:\; 0\le i\le 7\}$ (here the summation of subindexes are taken modulo 7), we show that  \begin{equation}\label{no-P4-cut}
\mbox{no element of } {\cal S} \mbox{ can be a cutset}.
\end{equation}

If (\ref{no-P4-cut}) does not hold, let $S=\{v_0, v_1, v_2, v_3\}$, and suppose that $S$ is a cutset. Let $G_1$ and $G_2$ be two subgraphs of $G$ with $V(G_1)\cap V(G_2)=S$ and $E(G_1)\cup E(G_2)=E(G)$, and suppose that $v_4\in V(G_2)$. Let $x$ be a vertex of $V(G_1)\setminus S$. By our assumption that every 3-cut of $G$ of size 3 must be a stable set, $G_1$ has four $(x, S)$-paths, say $P_0, P_1, P_2$ and $P_3$, such that $V(P_i)\cap S=v_i$ for $i\in\{0, 1, 2, 3\}$. By (\ref{theta+1}), $G_1$ has an induced $v_0v_2$-path, say $P$, of length at least 3 with $P^*\subseteq V(G_1)\setminus S$. But now, either $Pv_0v_7v_6v_2$ or $Pv_0v_4v_5v_6v_2$ is a big-odd-hole. This contradiction proves (\ref{no-P4-cut}).

\medskip

Let $S=\{v_0, v_1, v_2, v_3\}$,  $S'=\{v_4, v_5, v_6, v_7\}$, and let $Q$ be an induced $v_0v_2$-path other than $v_0v_1v_2$.  If $e(S', Q^*)=0$ then either $v_2v_6v_7v_0Q$ or $v_2v_6v_5v_4v_0Q$  is a big-odd-hole. Therefore,
\begin{equation}\label{v2-neighbors}
e(S', Q^*)\neq 0 \mbox{ for each induced $v_0v_2$-path } Q \mbox{ other than } v_0v_1v_2.
\end{equation}

By (\ref{no-P4-cut}), $S$ is not a cutset. Let $A_0$ be the set of all vertices on some $(v_1, v_3)$-bad ears, and let $A=A_0\setminus \{v_1, v_3\}$. Let $P$ be a shortest $(A, S')$-path in $G-S$. By Claim~\ref{bad-ears-1}, each vertex in $V(G)\setminus V(H')$ has at most one neighbor in $H'$, and so $\l(P)\ge 3$. Let $z$ be the end of $P$ in $S'$. Recall that $v_1w_1w_2v_3$ denotes a bad ear of $H$. Without loss of generality, we may suppose that $\{w_1, w_2\}\cap V(P)\neq \emptyset$. We will show that $P$ is a $(v_4, w_2)$-path.

\begin{claim}\label{not-w1}
$w_1\not\in V(P)$.
\end{claim}
\pf Suppose to its contrary that $w_1$ is an end of $P$.  Then, $e(P-z, S')=0$ by the shortestness of $P$. If both $v_0$ and $v_2$ have neighbors on $P$, then $P\cup\{v_0, v_2\}$ has an induced $v_0v_2$-path $P_0$ such that $e(S', P^*_0)=0$, contradicting (\ref{v2-neighbors}). So, we have that
\begin{equation}\label{v2-neighbors-1}
\mbox{either $v_0$ or $v_2$ has no neighbors on $P$}.
\end{equation}

First we suppose that $z=v_4$. If $v_2$ has neighbors on $P$, let $P_1$ be the shortest $v_2v_4$-path in $P\cup\{v_2\}$, then $v_0$ has no neighbors on $P_1$ by (\ref{v2-neighbors-1}), and so either $v_2v_6v_5v_4P_1$  or $v_2v_6v_7v_0v_4P_1$  is a big-odd-hole. Therefore, $v_2$ has no neighbors on $P$.  Let $P_2$ be the shortest $v_1z$-path in $P\cup\{v_1\}$, and let $y$ be the neighbor of $v_3$ on $P_2$ which is closest to $v_1$. If $y=v_4$ then either $v_1v_2v_3v_4P_2$  or $v_1v_2v_6v_5v_4P_2$  is a big-odd-hole. So, $y\neq v_4$. It is certain that $\l(P_2[v_1, y])\ge 3$ by the choice of $A$. By Claim~\ref{bad-ear-0},
$v_0$ has no neighbors in $P_2[v_1, y]-v_1$, and then either $v_1v_0v_4v_3yP_2[y, v_1]$  or $v_1v_2v_3yP_2[y, v_1]$  is a big-odd-hole. This shows that $z\ne v_4$.

\medskip

Next we suppose that  $z=v_5$. Let $P_1$ be the shortest $v_1v_5$-path in $P\cup \{v_1\}$. If $v_0$ has no neighbors on $P_1-v_1$,  then either $v_1v_0v_4v_5P_1$  or $v_1v_0v_7v_6v_5P_1$  is a big-odd-hole. So, $v_0$ has neighbors on $P_1-v_1$. Let $P_2$ be the shortest $v_0v_5$-path in $P_1\cup\{v_0\}$. Since $v_2$ has no neighbors on $P_2$ by (\ref{v2-neighbors-1}), we have that either $v_0v_7v_6v_5P_2$  or $v_0v_1v_2v_6v_5P_2$ is a big-odd-hole. Therefore, $z\ne v_5$.

\medskip

After that we suppose that $z=v_6$. Let $P_1$ be the shortest $v_3v_6$-path in $P\cup \{v_3, w_2\}$. If $v_0$ has no neighbors on $P_1$, then either $v_3v_4v_5v_6P_1$  or $v_3v_4v_0v_7v_6P_1$  is a big-odd-hole. Otherwise, $v_0$ has neighbors on $P_1$. Let $P_2$ be the shortest $v_0v_6$-path in $P_1\cup\{v_0\}$. By (\ref{v2-neighbors-1}), $v_2$ has no neighbors on $P_2$, and so either $v_0v_4v_5v_6P_2$  or $v_0v_4v_3v_2v_6P_2$ is a big-odd-hole. This shows that $z\ne v_6$.

\medskip

Finally, we suppose that $z=v_7$, and let $P_1$ be the shortest $v_3v_7$-path in $P\cup \{v_3, w_2\}$. If $v_2$ has no neighbors on $P_1$, then either $v_3v_2v_6v_7P_1$  or $v_3v_4v_5v_6v_7P_1$ is a big-odd-hole.
Otherwise, let $P_2$ be the shortest $v_2v_7$-path in $P_1\cup \{v_2\}$, then $v_0$ has no neighbors on $P_2$ by  (\ref{v2-neighbors-1}) which indicates that $v_1$ has no neighbors on $P_2$ as $H$ has no bad ears with attachments $\{v_1, v_7\}$. Now, either $v_2v_1v_0v_7P_2$  or $v_2v_3v_4v_0v_7P_2$  is a big-odd-hole. This contradiction proves Claim~\ref{not-w1}. \qed

\medskip

Next, we show that
\begin{claim}\label{v4-w2}
$P$ is a $v_4w_2$-path, and  $e(\{v_0, v_1, v_2, v_5, v_6, v_7, w_1, x_1, x_2\}, P^*)=0$.
\end{claim}
\pf From Claim~\ref{not-w1}, it is certain that $w_2$ must be an end of $P$. We show that $z\not\in \{v_5, v_6, v_7\}$.

Suppose that  $z=v_5$. It is certain that $N_{P^*}(v_i)=\emptyset$ for each $i\in\{4, 6, 7\}$.  If $N_{P}(v_2)=\emptyset$, then $N_{P^*}(v_3)=\emptyset$ by Claim~\ref{bad-ear-0} and by our assumption that $G$ has no $(v_3, v_5)$-bad ears, and so either $Pw_2v_3v_4v_5$ or $Pw_2v_3v_2v_6v_5$ is a big-odd-hole. Therefore, $N_{P}(v_2)\ne \emptyset$, and thus $N_{P}(v_0)=\emptyset$ by (\ref{v2-neighbors}). Let $x$ be the neighbor of $v_2$ on $P$ which is closest to $v_5$, and let $P'=P[x, v_5]$. Then, $N_{P'}(v_3)=\emptyset$ (again by Claim~\ref{bad-ear-0} and by our assumption that $G$ has no $(v_3, v_5)$-bad-ears). If $N_{P'}(v_1)=\emptyset$, then either $P'v_2v_3v_4v_5$ or $P'v_2v_1v_0v_4v_5$ is a big-odd-hole. Otherwise, let $x'$ be the neighbor of $v_1$ on $P'$ which is closest to $v_5$, and let $P''=P[x', v_5]$. Then, either $P''v_1v_2v_3v_4v_5$ or $P''v_1v_0v_4v_5$ is a big-odd-hole. This contradiction shows that $z\neq v_5$.

\medskip

Suppose that $z=v_6$. Let $P_1$ be the shortest $v_3v_6$-path in $P\cup \{v_3\}$. If $v_0$ has no neighbors on $P_1$, then either $v_3v_4v_5v_6P_1$  or $v_3v_4v_0v_7v_6P_1$  is a big-odd-hole. So, $v_0$ has neighbors on $P_1$. Otherwise, let $P_2$ be the shortest $v_0v_6$-path in $P_1\cup\{v_0\}$. Then, $v_2$ has no neighbors on $P_2-v_6$ by (\ref{v2-neighbors}), and so either $v_0v_4v_5v_6P_2$  or $v_0v_4v_3v_2v_6P_2$  is a big-odd-hole. Therefore, $z\ne v_6$.

\medskip

Suppose that $z=v_7$, and let $P_1$ be the shortest $v_3v_7$-path in $P\cup \{v_3, w_2\}$. If $v_2$ has no neighbors on $P_1$, then either $v_3v_2v_6v_7P_1$  or $v_3v_4v_5v_6v_7P_1$ is a big-odd-hole. Otherwise, let $P_2$ be the shortest $v_2v_7$-path in $P_1\cup \{v_2\}$. Then, $v_0$ has no neighbors on $P_2-v_7$ by  (\ref{v2-neighbors-1}), and so $v_1$ has no neighbors on $P_2-v_2$ as $H$ has no $(v_1, v_7)$-bad ears. Now, either $v_2v_1v_0v_7P_2$ or $v_2v_3v_4v_0v_7P_2$  is a big-odd-hole. Therefore, $z\ne v_7$.

\medskip

Now, we have that $z=v_4$, i.e., $P$ is a $(v_4, w_2)$-path. It is certain that $e(\{v_5, v_6, v_7, w_1\}, P^*)=0$. To prove Claim~\ref{v4-w2}, we need only check that $e(\{v_0, v_1, v_2, x_1, x_2\}, P^*)=0$.

If $v_2$ has neighbors on $P$, let $P_1$ be the shortest $v_2v_4$-path in $P\cup\{v_2\}$, then $v_0$ has no neighbors on $P_1$ by (\ref{v2-neighbors}), and so either $v_2v_6v_5v_4P_1$  or $v_2v_6v_7v_0v_4P_1$  is a big-odd-hole. Therefore, $v_2$ has no neighbors on $P$.

If $v_1$ has neighbors on $P$, let $P'_1$ be the shortest $v_1v_4$-path in $P\cup\{v_1\}$, then $v_3$ has no neighbors on $P'_1$ by the choice of $A$ and $P$, and so either $v_1v_2v_6v_5v_4P'_1$  or $v_1v_2v_3v_4P_1$ is a big-odd-hole. Therefore, $v_1$ has no neighbors on $P$.

If $v_0$ has neighbors on $P-v_4$, let $v'$ be the neighbor of $v_0$ on $P$ which is closest to $w_2$, and let $v''$ be the neighbor of $v_3$ on $P[w_2, v']$ which is closest to $v'$, then either $P[v', v'']v''v_3v_2v_1v_0v'$ or $P[v', v'']v''v_3v_2v_6v_7v_0v'$ is a big-odd-hole. This shows that $v_0$ has no neighbors on $P$ also.

If $x_1$ has neighbors on $P$, let $x'$ be the neighbor of $x_1$ on $P$ which is closest to $w_2$, and let $x''$ be the neighbor of $v_3$ on $P[w_2, x']$ which is closest to $x'$, then either $P[x', x'']x''v_3v_4v_5x_1x'$ or $P[x', x'']x''v_3v_2v_6v_5x_1x'$ is a big-odd-hole. With a similar argument, one can find a big-odd-hole if $x_2$ has neighbors on $P$. Therefore, Claim~\ref{v4-w2} holds. \qed

\medskip

Let $T=\{v_0, v_7, v_6, v_5\}$,  and let $T'=\{v_1, v_2, v_3, v_4\}$. Let $B_0$ be the set of all vertices on some $(v_5, v_7)$-bad ears, and let $B=B_0\setminus \{v_5, v_7\}$. Recall that $H$ has a bad ear $v_5x_1x_2v_7$. By applying the same arguments to $T, T'$, and $B$ as above used to $S, S'$ and $A$, we may suppose that $G$ has an induced $v_4x_1$-path $P'$ with $e(\{v_0, v_1, v_2, v_3, v_6, v_7, w_1, w_2, x_2\}, (P')^*)=0$. Let $Q$ be a shortest $w_2x_1$-path in $P\cup P'$. Combining with Claim~\ref{v4-w2}, we have either $Qw_2w_1v_1v_0v_7x_2x_1$ or $Qw_2w_1v_1v_2v_6v_7x_2x_1$ is a big-odd-hole. This contradiction completes the proof of Case 2. Consequently, we complete  the proof of Theorem~\ref{4-critical1}.  \qed

\section{Remarks}

Recall that ${\cal G}_{\l}$ denotes the family of graphs that have girth $2\l+1$ and have no odd holes of length at least $2\l+3$, $\theta^+$ denotes the graph obtained from the Petersen graph by removing two adjacent vertices, $\theta$ denotes the graph obtained from $\theta^+$ by removing an edge incident with two 3-vertices, and $\theta^-$ is the graph obtained from $\theta^+$ by deleting a 2-vertex.

Theorem~\ref{coro-4-colorable} asserts that every graph of ${\cal G}_2$ inducing no 5-cycles sharing edges is 3-colorable. If $G\in {\cal G}_2$ induces 5-cycles sharing edges, then $G$ must induce $\theta$ or $\theta^-$. To confirm Conjecture~\ref{conj-1}, it suffices now to prove Conjecture~\ref{conj-SWXX}, i.e., graphs of ${\cal G}_2$ inducing neither $\theta$ nor $\theta^-$ must be 3-colorable.

Let $G$ be a graph  in ${\cal G}_{\l}$, $\l\ge 2$. Theorem~\ref{lem-subgraph} asserts that $\chi(G)\le 4$, and that for each vertex $u$ of $G$, all the vertices of distance $i$ ($i\ge 0$) from $u$ induce a bipartite subgraph. Assuming Conjecture~\ref{conj-1}, it is reasonable to believe that $\chi(G)\le 3$ for each $G\in{\cal G}_{\l}$, where $\l\ge 2$.
\begin{conjecture}\label{conj-2}
$\chi(G)\le 3$ for each graph $G\in \cup_{\l\ge 2}{\cal G}_l$.
\end{conjecture}

\bigskip

\noindent{\bf Acknowledgement}: We thank Dr. Jie Ma for pointing out the 3-connectivity of graphs in ${\cal G}_0$, and thank Dr. Sophie Spirkl for her helps on the proof of Theorem~\ref{coro-4-colorable} and on the presentation.

\end{document}